\documentclass[12pt,notitlepage]{article}[1995/06/26]
\usepackage{a4wide,latexsym,amsfonts,amssymb,enumerate,amsmath}
\usepackage[dvips]{graphics}
\usepackage{natbib}
\usepackage[english]{babel}

\newtheorem{defn}{Definition}[section]
\newtheorem{lem}[defn]{Lemma}
\newtheorem{thm}[defn]{Theorem}
\newtheorem{cor}[defn]{Corollary}
\newtheorem{assu}[defn]{Assumption}
\newcommand{\cip}{\mbox{$\perp\!\!\!\perp$}}
\newcommand{\proof}{\bigskip\noindent{\bf Proof.\enskip}}
\newcommand{\proofof}[1]{\bigskip\noindent{\bf Proof of #1.\enskip}}
\newcommand{\Endproof}{\ \hfill$\Box$\bigskip}
\newcommand{\GG}{g}

\begin{document}

\title{Estimating the causal effect of a time-varying treatment on
time-to-event using structural nested failure time models}

\author{Judith Lok, Richard Gill, Aad van der Vaart and James Robins\\
University of Leiden, Utrecht University, Vrije Universiteit Amsterdam\\
and Harvard University}

\date{July 2003}

\maketitle

\begin{abstract}
In this paper we review an approach to estimating the causal
effect of a time-varying treatment on time
to some event of interest. This approach is designed for
the situation where the treatment may have been repeatedly adapted to
patient characteristics, which themselves may
also be time-dependent. In this situation the effect of
the treatment cannot simply be estimated by conditioning on the patient
characteristics, as these may themselves be indicators of the
treatment effect. This so-called time-dependent confounding is typical
in observational studies. We discuss a new class of failure
time models, structural nested failure time models, which can be used
to estimate the causal effect of a time-varying treatment,
and present methods for estimating and testing the parameters of these models. 
\end{abstract}

\section{Introduction}
This paper offers a new approach to estimating, from observational
data, the causal effect of a time-dependent treatment on time to an
event of interest in the presence of time-dependent confounding
variables. This approach is based on a new class of failure time
models, the \emph{structural nested failure time models} (SNFTM). The
primary goal of this paper is to motivate the need for structural
nested failure time models. To achieve this goal in the most
straightforward manner, we shall assume that the event times are
observed without censoring, and that there is no missing or
misclassified data. Additional complications that arise when these
assumptions are not satisfied are discussed in Robins et al.~(1992)
and Robins~(1993).

The approach using SNTFMs will be useful in any observational study in
which there exist time-dependent risk factors that are also predictive
for subsequent exposure to the treatment under study, i.e.\ in any
study where there are time-dependent covariates that correlate with
the final outcome of the treatment, but also with the amount or type
of treatment over time.  This situation arises in any observational
study in which there is ``treatment by indication'', i.e.\ the
treatment is not predetermined by the investigator, but adapted to the
current condition of the patient.  The problem then is to distinguish
between treatment effect and selection bias (i.e.\ confounding). For
example, in an observational study for the effect of AZT treatment on
HIV-infected subjects, subjects with low CD4 lymphocyte counts at a
given time are subsequently at increased risk of developing AIDS and
are for that reason more likely to be treated with AZT. Thus the
covariate variables ``low CD4-count'' is a risk factor for AIDS, but
is also a predictor of subsequent treatment with AZT. The problem is
then to isolate the effect of AZT treatment as given according to a
predetermined plan (which may take into account covariates) from the
confounding effect of CD4-count.  As a second example, many physicians
withdraw women from exogenous estrogens at the time they develop an
elevated blood cholesterol, since both exogenous estrogens and
elevated blood cholesterol are considered possible cardiac risk
factors. Therefore, in a study of the effect of postmenopausal
estrogen on cardiac mortality, the covariate variables ``cholesterol
level'' is a predictor of subsequent exposure to estrogens, but also
correlates with the outcome ``cardiac mortality''. As a third example,
in observational studies of the efficacy of cervical cancer screening
on mortality, women who have had operative removal of their cervix due
to invasive disease are no longer at risk for further screening (i.e.\
exposure), but are at increased risk for death. Therefore, the
covariate, ``operative removal of the cervix'', is an independent risk
factor for death, but also a predictor of subsequent exposure. As a
final epidemiologic example, in occupational mortality studies,
unhealthy workers who terminate employment early are at increased risk
of death compared to other workers and receive no further exposure to
the chemical agent under study. Therefore, the time-dependent
covariate ``employment status'' is an independent risk factor for
death, and a predictor of exposure to the study agent.

Epidemiologists refer to the covariates in the preceding
examples as ``time-dependent confounders''. It may be
important to analyze the data from any of the above studies using the
approach presented in this paper.

For pedagogic purposes, we shall illustrate our models and assumptions
throughout the paper by the problem of estimating, from data obtained
in an observational study, the effect of treatment with the drug AZT
on time to clinical AIDS in asymptomatic subjects with newly diagnosed
human immunodeficiency virus (HIV) infection. We shall suppose that
measurements on current AZT dosage as well as on various
time-dependent covariates, such as weight, temperature, hematocrit,
and CD4-lymphocyte count, are recorded at regularly spaced time points,
until the development of clinical AIDS. These time points, which we
denote by $0=\tau_0<\tau_1<\tau_2<\cdots<\cdots<\tau_K$, may for
instance correspond to clinic visits at which the measurements are
obtained, with time defined as time since the diagnosis of HIV
infection.

Our goal will be to identify and estimate, for each \emph{treatment
regime}, the time-to-AIDS distribution that would have been observed
if (typically \emph{counter to fact}) each study subject had followed
the AZT treatment history prescribed by the regime. We shall call each
such distribution an AZT treatment regime-specific, counterfactual,
time-to-AIDS distribution. The treatment regimes we study need not be
static. A \emph{treatment regime} is a rule that assigns to each
possible covariate history through time $\tau_k$, an AZT dosage rate
$a_k$ to be taken in the interval $\left(\tau_k,\tau_{k+1}\right]$. A
simple example of a treatment regime is ``take an AZT dosage $a_k$ of
$1,000$ milligrams of AZT daily in the interval
$\left(\tau_k,\tau_{k+1}\right]$ if the hematocrit measured at
$\tau_k$ exceeds $30$; otherwise take no AZT in the interval''.

Our interest in AZT treatment regime-specific, counterfactual
time-to-AIDS distributions is based on the following
considerations. Suppose, after the completion of the study, a further
individual with newly diagnosed HIV infection, whom we shall call
``the infected subject'', wishes to use the data from the completed
study to select the AZT dosage schedule that will maximize his
expected or median number of years of AIDS-free survival. If the
``infected subject'' is considered exchangeable with the subjects in
the trial, then he would wish to follow the AZT treatment regime whose
regime-specific, counterfactual time-to-AIDS distribution has the
largest expected or median value.

In Section~\ref{Gcomps} we show that the AZT treatment
regime-specific, counterfactual time-to-AIDS distributions are
identified from the observed data under the assumption that the
investigator has succeeded in recording sufficient data on the history
of all covariates to ensure that, at each time $\tau_k$, given the
covariate history and the AZT treatment history up till $\tau_k$, the
AZT dosage rate in $\left(\tau_k,\tau_{k+1}\right]$ is independent of
the regime-specific, counterfactual time-to-AIDS. Robins~(1992)
refers to this assumption as the assumption of \emph{no unmeasured
confounding factors}. In other words, under this assumption at each
time point the treatment can be viewed as depending only on recorded
information up till that point and external factors that are not
predictive of (counterfactual) survival.

In Section~\ref{repars} we introduce \emph{structural nested failure
time models (SNFTM)}. An SNFTM models the magnitude of the causal
effect of a (final) blip of AZT treatment in the interval
$\left(\tau_k,\tau_{k+1}\right]$ on time-to-AIDS, as a function of
past AZT and covariate history.  We show that, under the assumption of
no unmeasured confounding, the null hypothesis of no causal effect of
AZT on time-to-AIDS is equivalent to the null hypothesis that the
parameter vector of any SNFTM is $0$. 

The term ``structural'' in SNFTM derives terminology used
in the social science and econometric literature (e.g.\ Rubin~(1978)).
Our models are ``structural'', because they
directly model regime-specific, counterfactual time-to-AIDS
distributions. In Sections~\ref{mlesec} and~\ref{Gest} we discuss
two different methods to fit SNFTMs and to use them for inference.

In Section~\ref{mlesec} we show that, under the assumption of no
unmeasured confounding, SNFTMs can be understood as a component of a
particular reparameterization of the joint distribution of the
observables. We use this reparameterization to develop
likelihood-based tests of the causal null hypothesis of no effect of
AZT-exposure on time-to-AIDS. We also show how to estimate the
AZT-treatment regime-specific, counterfactual time-to-AIDS
distributions, in the case that the null hypothesis of no causal
effect of AZT on time-to-AIDS is rejected.

In Section~\ref{Gest} we present an alternative, semiparametric
approach to test the null hypothesis of no treatment effect and to
estimate the parameters in an SNFTM. This approach,
\emph{G--estimation}, has the advantage of avoiding for
parameterization of the distributions appearing in the
likelihood-based approach of Section~\ref{mlesec} (e.g. the
conditional distributions of covariates given past treatment- and
covariate history). Instead G--estimation uses a model
for the SNFTM and for the conditional distribution of treatment given
past treatment- and covariate history. Tests and estimators based on
G-estimation have the additional advantage that they can often be
calculated with standard software.

\section{Formalization of the problem}
We fix a discrete time frame $\tau_0=0<\tau_1<\tau_2<\ldots<\tau_K$
throughout the paper, where $\tau_0$ is the time of enrollment in the
study (and possibly also initiation of treatment), $\tau_1,\tau_2,\ldots$ are
the times of the clinic visits, and $\tau_K$ can be the time of the
last clinic visit, or can be chosen past the upper support point of
the time-to-AIDS distribution. For simplicity the times of the clinic
visits are assumed to be the same for all patients (as long as they
are alive).

At each time point $\tau_k$ we measure a covariate vector $L_k$ for
each patient, where $L_0$ may also contain time-independent covariates
and information collected before time $\tau_0$, and we register the
treatment given in the interval $(\tau_k, \tau_{k+1}]$ in a variable
$A_k$, for instance the AZT dosage, assumed constant during the
interval. Besides covariates $L_k$ and treatments $A_k$, we observe
for each person a positive time $T$, for instance the time from
enrollment to the development of clinical AIDS. Thus the data observed
on one person is a vector $(\overline L_K,\overline A_K,T)$, where,
for each $k=0,1,\ldots,K$,
\begin{eqnarray*}
\overline L_k&=&(L_0,L_1,\ldots, L_k),\\
\overline A_k&=&(A_0,A_1,\ldots, A_k).
\end{eqnarray*}
For time instances $\tau_k>T$ the values $L_k$ and $A_k$ may be
interpreted to be empty. For simplicity we assume that the variables
$L_k$ and $A_k$ take their values in countable sets, denoted by
${\mathcal L}_k$ and ${\mathcal A}_k$. The total set of observations
are a sample of $n$ independent and identically distributed (i.i.d.)
observations from the distribution of the random vector $(\overline
L_K,\overline A_K,T)$.

As is clear from the preceding display we use the overline notation
$\overline{}$ to denote a ``cumulative vector''. For simplicity of
notation, it will be understood that whenever two expressions such as
$\overline l_k$ and $\overline l_{k-1}$ occur together, then
$\overline l_{k-1}$ is the initial part of $\overline l_k$.

A ``treatment regime'' is a prescription for the treatment dosages
fixed at the times $\tau_k$, where at each time instant the prescribed
treatment may depend on the observed covariate history until this
time. We make this precise in the following definition.

\begin{defn} \emph{(treatment regimes).}
A treatment regime $\GG$ is a vector  $\GG=(\GG_0,\ldots, \GG_K)$
of functions 
$\GG_k: {\mathcal L}_0\times\cdots\times{\mathcal L}_k\to {\mathcal A}_k$.
\end{defn}

The value $a_k=\GG_k(\overline l_k)$ of the $k$th coordinate 
of the treatment regime $\GG$ at covariate $\overline l_k$
is interpreted as the dosage
prescribed by treatment regime $\GG$ in the interval $(\tau_k,\tau_{k+1}]$
to a patient with covariate history $\overline l_k$ following 
this regime (up to time $\tau_k$). The treatment at time $\tau_k$
may depend on the full covariate history $\overline l_k=(l_0,\ldots,l_k)$
until time $\tau_k$, not just on $l_k$.
We define maps 
$\overline g_k: {\mathcal L}_0\times\cdots\times{\mathcal L}_k\to 
{\mathcal A}_0\times\cdots\times{\mathcal A}_k$ by
$$\overline g_k(\overline l_k) =\bigl(g_0(l_0),g_1(\overline
l_1),\ldots, g_k(\overline l_k)\bigr).$$ 
To alleviate notation we may drop the subscripts $k$ or the overline
in $g_k$ or $\overline g_k$ if the value of $k$ is clear from the
context.  In particular $g(\overline l_K)=\overline g(l_K)= \overline
g_K(\overline l_K)$ are equivalent notations for the complete
treatment history.

We wish to study the effect of treatment using the observed data.
Depending on this data not all treatment regimes may be accessible to
analysis. We call a treatment regime ``evaluable'' (relative to the
distribution of the data vector $(\overline L_K,\overline A_K,T)$) if
whenever the regime was followed until some time $\tau_k$ by some
positive fraction of the population, then it is also followed in the
interval $(\tau_k, \tau_{k+1}]$.

\begin{defn} \emph{(evaluable treatment regimes).}
A treatment regime $\GG$ is called evaluable if for each $k$ and
each $\overline{l}_k\in\overline{{\mathcal L}}_k$,
\begin{equation*}
P\left(\overline{L}_k=\overline{l}_k,\overline{A}_{k-1}=
\overline{\GG}\left(\overline{l}_{k-1}\right),T>\tau_k\right)>0 \Rightarrow
P\left(\overline{L}_k=\overline{l}_k,\overline{A}_k=
\overline{\GG}\left(\overline{l}_k\right),T>\tau_k\right)>0.
\end{equation*}
\end{defn}

Next we introduce \emph{counterfactual variables}. These will be
instrumental both to express the aims of the statistical analysis, and
to formulate our assumptions. In our mathematical model the
counterfactual variables are ordinary random variables $T^\GG$,
one for each treatment regime $\GG$, that are assumed to be defined on the
same probability space as the data vector $(\overline L_K,\overline
A_K,T)$.  The variable $T^{\GG}$ should be thought of as a patient's
time to clinical AIDS had she been treated according to treatment
regime $\GG$. Because in actual fact the patient receives treatment
$\bar A_K$ (resulting in time to aids $T$), the variable $T^\GG$ is ``counter
to fact''. However, it gives a useful notation to express the
distribution of interest, and will be related to the observable
variables by two assumptions.

Counterfactual variables referring to different subjects are
assumed independent (cf.\ Rubin~(1978)), and hence we can formulate
our set-up in terms of the set of random variables $(T^\GG, T, \overline
L_K,\overline A_K)$ referring to one person. We shall not be
interested in the joint distribution of counterfactual variables
corresponding to different treatment regimes.  We also do not need
counterfactual versions of the covariates or treatments.

We describe the aims of the statistical analysis in terms of the
counterfactual variables.  The \emph{G--null hypothesis} of no effect
of AZT on time-to-AIDS is the hypothesis that
\begin{equation*}
P\left(T^{g_1}>t\right)=P\left(T^{g_2}>t\right)\hspace{0.7cm}
{\rm for}\;{\rm all}\; {\rm treatment}\;{\rm regimes}\;g_1\;{\rm and}\;g_2.
\label{nte}
\end{equation*}
In Section~\ref{mlesec} we derive fully parametric
likelihood-based tests of this G--null hypothesis based on  a random
sample from the distribution of the
observables $\left(\overline{L}_K,\overline{A}_K,T\right)$,
and a parametric model for their joint
distribution. In Section~\ref{Gest} we develop an alternative,
semi-parametric procedure with the same aim.

If the G--null hypothesis is rejected, then the next goal is to
identify and estimate, for each treatment regime $g$, the survival
curve $t\mapsto P\left(T^\GG>t\right)$, i.e.\ the survival curve that
would have been observed had a subject followed regime
$g$. Specifically, if our infected subject outside of the study
mentioned in the introduction wishes to maximize his expected years of
AIDS-free survival, he would follow the regime $g$ that maximized $E
T^\GG=\int_0^\infty P\left(T^\GG>t\right)\, dt$. Inference regarding
the distribution of counterfactual variables is referred to as
\emph{causal inference}, as the outcomes $T^\GG$ are interpreted
as being the effect of the treatment regime $\GG$.

Clearly it is impossible to make inference about the
counterfactual survival distributions $P(T^\GG>t)$ based on the
observed data unless the variables $T^\GG$ and 
$(\overline L_K,\overline A_K,T)$
are related.  The assumed coupling of these variables
on a given underlying probability space allows to make
the following assumptions relating counterfactual and factual
variables.

\begin{assu}\emph{(consistency).} \label{cons}
For any treatment
regime $\GG$, $\overline{l}_k\in\overline{{\mathcal L}}_k$ and
$t\in\left(\tau_k,\tau_{k+1}\right]$,
\begin{equation*}
\left\{T^{\GG}>t,\overline{L}_k=\overline{l}_k,\overline{A}_k=
\overline{\GG}\left(\overline{l}_k\right),T> \tau_k\right\} =
\left\{T>t,\overline{L}_k=\overline{l}_k,\overline{A}_k=
\overline{\GG}\left(\overline{l}_k\right),T> \tau_k\right\}.
\end{equation*}
\end{assu}
\begin{assu}\emph{(no unmeasured confounding).} \label{ass}
For any treatment regime $\GG$, for any time $\tau_k$ and for
any $\overline{l}_k\in \overline{{\mathcal L}}_k$,
\begin{equation*} A_k \cip T^{\GG} | \overline{L}_k=\overline{l}_k,
\overline{A}_{k-1} = \overline{\GG}\left(\overline{l}_{k-1}\right).
\end{equation*}
\end{assu}

Here the notation $X\cip Y|Z=z$, borrowed from Dawid~(1979),
means that the random variabless $X$ and $Y$ are conditionally
independent given the event $Z=z$.

The consistency assumption, Assumption~\ref{cons},  couples the true and
counterfactual survival times $T$ and $T^{\GG}$ by merely stating that
if until some time $\tau_k$ a patient is treated exactly as prescribed
by regime $g$, then she would die at some time in the interval
$(\tau_k,\tau_{k+1}]$ under regime $g$ if and only if she actually
died at the same time. This implies in particular that if all patients
were treated according to a predetermined treatment regime, then
counterfactual and actual survival times coincide. This is the customary
situation in clinical trials, but may fail to be the case
in an observational study.

The assumption of no unmeasured confounding, Assumption~\ref{ass}, can
be expected to hold if the observed covariate history $\overline L_K$
contains sufficient information, so that at each time $\tau_k$ the
treatment $A_k$ can be assumed to depend on the covariate history
$\overline L_k$ of a patient up till that time and no other relevant
information. The assumption would for instance hold if at each time
$\tau_k$ the treatment in the interval $(\tau_k,\tau_{k+1}]$ is
assigned through randomization within fixed levels of equal covariates
$\overline L_k$ and earlier treatments.

More specifically, in our AIDS example Assumption~\ref{ass} may be
expected to hold if the following information is recorded in
$\overline{L}_k$: all risk factors (i.e.\ predictors) of
regime-specific, counterfactual time-to-AIDS, other than prior
AZT-history $\overline{A}_{k-1}$, that are used by physicians and
patients to determine the dose $A_k$ of AZT in
$\left(\tau_k,\tau_{k+1}\right]$. Then, given $\overline{L}_k$ and
$\overline{A}_{k-1}=\GG\bigl(\overline{L}_{k-1}\bigr)$, the treatment
$A_k$ in the interval $(\tau_k,\tau_{k+1}]$ may
be thought of as depending only on external factors unrelated to
the patient's prognosis regarding time-to-AIDS, and hence as being
independent of $T^\GG$. For example, since it is known that physicians
tend to prescribe AZT to subjects with low CD4-counts and a low
CD4-count is an independent predictor of time-to-AIDS, the assumption
of no unmeasured confounding would be false if $\overline{L}_k$ does
not contain CD4-count history. 

It is a basic objective of epidemiologists conducting an observational
study to collect data on a sufficient number of covariates to ensure
that Assumption~\ref{ass} will be true. In this paper, we assume this
objective has been realized, while recognizing that, in practice, this
may only approximately be the case.

\section{G--computation}\label{Gcomps}
We are interested in the distribution of the counterfactual, and hence
unobservable, variables $T^{\GG}$, as they indicate the success or
failure from applying the treatment regime $\GG$.  In this section we
show that, under Assumptions~\ref{cons} and~\ref{ass}, the
distribution of $T^\GG$ is identifiable from the distribution of the
observed data $\bigl(\overline{L}_K,\overline{A}_K,T\bigr)$ for each
evaluable treatment regime $\GG$.  As a consequence, given a random
sample from the latter distribution, the distribution of $T^\GG$ is
estimable, in principle.

In fact, the following \emph{G--computation formula} gives an explicit
expression for $P\left(T^{\GG}>t\right)$, as well as several
conditional survival functions, in terms of the distribution of the
data $\left(\overline{L}_K,\overline{A}_K,T\right)$.

\begin{thm}\label{Gcomp}
\emph{(G--computation-formula).}
Suppose that Assumptions~\ref{ass} (no
unmeasured confounding) and~\ref{cons} (consistency) hold, and that
$\GG$ is an evaluable treatment regime. Then for any $t>0$,
with $p$ defined by $\tau_p<t\le \tau_{p+1}$,
\begin{eqnarray}\label{Gf}
P\left(T^{\GG}>t\right) &=& \sum_{l_{0}} \cdots\sum_{l_{p-1}}\sum_{l_p}
\Bigg[P\Bigl(T>t| \overline{L}_{p}=\overline{l}_{p},
\overline{A}_{p}=\overline{\GG}\left(\overline{l}_{p}\right),
T>\tau_{p}\Bigr)
\nonumber\\
&&\hspace{2.0cm}\times\prod_{m=0}^{p}\Big\{
P\Bigl(T>\tau_m|\overline{L}_{m-1}=\overline{l}_{m-1},
\overline{A}_{m-1}=\overline{\GG}\left(\overline{l}_{m-1}\right),
T>\tau_{m-1}\Bigr)\nonumber\\
&&\hspace{2.2cm}
\times P\Bigl(L_m=l_m|
\overline{L}_{m-1}=\overline{l}_{m-1},
\overline{A}_{m-1}=\overline{\GG}\left(\overline{l}_{m-1}\right),
T>\tau_{m}\Bigr) \Big\}\Bigg].\nonumber
\end{eqnarray}
\end{thm}

In the preceding theorem we interpret variables indexed by $-1$ as not
present, and events concerning only such variables as being empty. For
instance, the conditional probability 
$P\bigl(L_m=l_m| \overline{L}_{m-1}=\overline{l}_{m-1},
\overline{A}_{m-1}=\overline{\GG}\left(\overline{l}_{m-1}\right),
T>\tau_{m}\bigr)$ is to be read as the probability $P(L_0=l_0)$ when
$m=0$.

All conditional probabilities on the right side concern
observable variables. Hence the theorem gives an explicit
description of the survival function of the counterfactual
variable $T^\GG$ in terms of the distribution of
the data $(\overline L_K, \overline A_K, T)$. 

It is instructive to evaluate the formula in the simple case that
$K=1$, when there exists only one treatment $A_0$ applied in the
single interval $(0,\tau_1]$. Then the G--computation formula yields,
for $t>0$,
$$P(T^\GG>t)=\sum_{l_0} P\bigl(T>t|L_0=l_0,A_0=g(l_0)\bigr)
\,P(L_0=l_0).$$
This shows that in general the distribution of the counterfactual
variable $T^\GG$ differs from the distribution of $T$, which can be
written in the form
$$P(T>t)=\sum_{l_0} P\bigl(T>t|L_0=l_0\bigr)\,P(L_0=l_0).$$ 
This difference is not too surprising, because the variable $T^\GG$
refers to the treatment regime $\GG$, whereas $T$ relates to the
observed outcomes under the actual treatments. Had all patients
received treatment $g$, then the two distributions would coincide.
More notable is the difference between the conditional distribution of
$T$ given $A_0=a_0$ and the distribution of $T^\GG$ for the fixed
treatment regime $\GG$ that assigns all patients to treatment $a_0$,
i.e.\ $g(l_0)=a_0$.  These two survival distributions can be written
\begin{eqnarray*}
P(T^{a_0}>t)&=&\sum_{l_0} P\bigl(T>t|L_0=l_0,A_0=a_0\bigr)\,P(L_0=l_0),\\
P\bigl(T>t| A_0=a_0\bigr)&=&\sum_{l_0} P\bigl(T>t|L_0=l_0,A_0=a_0\bigr)
\,P\bigl(L_0=l_0| A_0=a_0\bigr).
\end{eqnarray*}
The conditional distribution of $T$ given $A_0=a_0$ is estimable, in
principle, by taking only those patients into account who happened to
receive treatment $a_0$. The outcome distribution of this subset of
patients may however be different from the distribution of the
counterfactual variable $T^{a_0}$, as a result of ``selection bias''.
In the actual world some patients may be assigned other treatments
than $a_0$, where the assignment $A_0$ may correlate with the
covariate variable $L_0$. Therefore, the conditional and unconditional
distributions of $L_0$ given $A_0$ may differ, and consequently so may
the right hand sides of the display. It is the counterfactual survival
function $t\mapsto P(T^{a_0}>t)$ that is the relevant one to judge the
causal effect of treatment $a_0$. Randomization of treatment over
patients within fixed levels of the covariate would have made $L_0$
and $A_0$ independent, and the difference would disappear.  The
protocol of a controlled experiment may include such randomization,
but in a observational study it cannot be taken for granted. The
G--computation formula then shows, under some assumptions, how we can
still compute the relevant outcome distributions from the observed
data distribution.

We can make further comparisons after deriving
a similar representation for conditional probabilities involving
the counterfactual variables.

\begin{thm}\label{Gcompcond}
\emph{(G--computation-formula).}
Under the assumptions of Theorem~\ref{ass},
for any $k\in \{0,1,2,\ldots,K\}$ and any $\overline l_k$ such that
$P\left(\overline{L}_k=\overline{l}_k,
\overline{A}_{k-1}=\overline{\GG}\left(\overline{l}_{k-1}\right),
T>\tau_k\right)>0$,
for any $t>\tau_k$, and with $p\ge k$ defined by $\tau_p<t\le \tau_{p+1}$,
\begin{eqnarray}\label{Gfc}
\lefteqn{P\Bigl(T^{\GG}>t|\overline{L}_k=\overline{l}_k,
\overline{A}_{k-1}=\overline{\GG}\left(\overline{l}_{k-1}\right),
T>\tau_k\Bigr)}\nonumber\\ &=&\sum_{l_{k+1}} \cdots\sum_{l_{p-1}}\sum_{l_p}
\Bigg[P\Bigl(T>t|\overline{L}_{p}=\overline{l}_{p},
\overline{A}_{p}=\overline{\GG}\left(\overline{l}_{p}\right),
T>\tau_{p}\Bigr)\nonumber\\
&&\hspace{1.9cm}\times\prod_{m=k+1}^{p}\Big\{
P\Bigl(T>\tau_m|\overline{L}_{m-1}=\overline{l}_{m-1},
\overline{A}_{m-1}=\overline{\GG}\left(\overline{l}_{m-1}\right),
T>\tau_{m-1}\Bigr)\nonumber\\ &&\hspace{1.9cm}
\times P\Bigl(L_m=l_m|
\overline{L}_{m-1}=\overline{l}_{m-1},
\overline{A}_{m-1}=\overline{\GG}\left(\overline{l}_{m-1}\right),
T>\tau_{m}\Bigr) \Big\}\Bigg].
\end{eqnarray}
\end{thm}

Again variables indexed by $-1$ should be read as not being present.
Furthermore, a repeated summation of the form $\sum_{l_{k+1}}\cdots
\sum_{l_p}a_{k,p}(\overline l_{k},l_{k+1},\ldots,l_p)$ is considered
to be the single term $a_{k,k}(\overline l_k)$ if $k=p$, whereas the
product $\prod_{k+1}^p$ is to be read as 1 in this case. The summation
may be restricted to terms whose conditioning events have positive
probability.

Again we may evaluate this formula in the simple case 
of a single treatment interval. Then the formula in the
preceding theorem (with $k=0=p, K=1$) reduces to 
$$P\bigl(T^\GG>t| L_0=l_0\bigr) =
P\bigl(T>t|L_0=l_0,A_0=g(l_0)\bigr).$$ 
The right side is precisely the conditional distribution of the actual
survival time for a subject with covariate $l_0$ following the
treatment regime $g$. Intuitively, the conditional probabilities
$P\bigl(T>t|L_0=l_0,A_0=g(l_0)\bigr)$ are the correct ones for
evaluating the quality of treatment $g$ for a subject with covariate
value $l_0$, and the equality in the preceding display is actually a
direct consequence of the Assumptions~\ref{cons} and~\ref{ass}
relating the counterfactual and factual survival times. (We may add
$A_0=g(l_0)$ in the conditioning event on the left by
Assumption~\ref{ass}, and next use Assumption~\ref{cons} to see that
$T^\GG$ may be replaced by $T$.)

Henceforth, we shall denote the right side of (\ref{Gfc}) by
$s_{\overline{l}_k,\GG}\left(t\right)$. For $k=-1$ this reduces to the
right side in Theorem~\ref{Gcomp}, and we write it as
$s_{\GG}(t)$, interpreting $\overline{l}_{-1}$ as empty.  
Then Theorems~\ref{Gcomp}-\ref{Gcompcond}
can be reformulated as saying that under
Assumptions~\ref{cons} (consistency) and~\ref{ass} (no unmeasured
confounding), for every evaluable treatment regime $\GG$,
\begin{equation*}P\left(T^{\GG}>t\right)=s_{\GG}(t)
\end{equation*}
and, for every $k=0,1,\ldots, K$,
\begin{equation*}
P\Bigl(T^{\GG}>t|\overline{L}_k=\overline{l}_k,
\overline{A}_{k-1}=\overline{\GG}\left(\overline{l}_{k-1}\right),
T>\tau_k\Bigr)=s_{\overline{l}_k,\GG}(t).
\end{equation*}
These functions are survival functions of distributions
that concentrate on $(\tau_k,\infty)$.

Inspection of the G--computation formula shows that
$s_{\overline{l}_k,\GG}$ is a (complicated) function of the
distribution of the data vector
$\left(\overline{L}_K,\overline{A}_K,T\right)$ and depends on this
distribution only through the conditional distributions of the
covariates and the survival time given the past, given by
\begin{equation}\label{plm}
P\left(L_m=l_m|\overline{L}_{m-1}=\overline{l}_{m-1},
\overline{A}_{m-1}=\overline{a}_{m-1},T>\tau_m\right),
\end{equation}
and
\begin{equation}\label{ptm}
P\left(T>t|\overline{L}_{m-1}=\overline{l}_{m-1},
\overline{A}_{m-1}=\overline{a}_{m-1},
T>\tau_{m-1}\right).
\end{equation}
In particular, the functions $s_{g,\overline{l}_k}$ do not depend on
conditional laws of the treatment variables $A_m$ given the past.

\proofof{Theorems~\ref{Gcomp} and~\ref{Gcompcond}} 
We prove Theorems~\ref{Gcomp} and~\ref{Gcompcond} by backward
induction on $k$, for fixed $t$ (and hence also fixed $p$).  Formula
(\ref{Gfc}) with $k=-1$ can be read as the formula given by
Theorem~\ref{Gcomp}, so we restrict to proving (\ref{Gfc}).

For $k=p$ the left side of (\ref{Gfc}) is equal to
\begin{eqnarray*}
\lefteqn{P\left(T^{\GG}>t|\overline{L}_p=\overline{l}_p, \overline{A}_{p-1}
=\overline{\GG}\left(\overline{l}_{p-1}\right),T>\tau_p\right)}\\
&&\hspace{2.5cm}=
P\left(T^{\GG}>t|\overline{L}_p=\overline{l}_p, \overline{A}_{p}
=\overline{\GG}\left(\overline{l}_{p}\right),T>\tau_p\right)\\
&&\hspace{2.5cm}=
P\left(T>t|\overline{L}_p=\overline{l}_p, \overline{A}_{p}
=\overline{\GG}\left(\overline{l}_{p}\right),T>\tau_p\right),
\end{eqnarray*}
where in the first equality we can add
$A_p=\GG_p\left(\overline{l}_p\right)$ in
the conditioning event by Assumption~\ref{ass} of no unmeasured confounding,
and in the second equality we can replace the event $T^{\GG}>t$
by the event $T>t$, because of the Assumption~\ref{cons}
of consistency.

The induction step is proved by similar arguments. Supposing
that (\ref{Gfc}) holds for $k\le p$,
we shall deduce that it also holds for $k-1$. We have
\begin{eqnarray*}
\lefteqn{P\left(T^{\GG}>t|\overline{L}_{k-1}=\overline{l}_{k-1},
\overline{A}_{k-2}=\overline{\GG}\left(\overline{l}_{k-2}\right),
T>\tau_{k-1}\right)}\\
&=&
P\left(T^{\GG}>t|\overline{L}_{k-1}=\overline{l}_{k-1},\overline{A}_{k-1}=
\overline{\GG}\left(\overline{l}_{k-1}\right),T>\tau_{k-1}\right)\\
&=&
P\left(T^{\GG}>\tau_{k}|\overline{L}_{k-1}=\overline{l}_{k-1},
\overline{A}_{k-1}=\overline{\GG}\left(\overline{l}_{k-1}\right),
T>\tau_{k-1}\right)\\
&&\hspace{1.5cm} \times P\left(T^{\GG}>t|\overline{L}_{k-1}=\overline{l}_{k-1},
\overline{A}_{k-1}=\overline{\GG}\left(\overline{l}_{k-1}\right),
T>\tau_{k-1},T^{\GG}>\tau_{k}\right).
\end{eqnarray*}
The first equality follows by the assumption of no unmeasured
confounding, while the second follows by conditioning on the
event $T^{\GG}>\tau_{k}$, where we note that $t>\tau_{k}$, because
$t>\tau_p\ge\tau_{k}$.  By the consistency assumption we can replace
the event $T^{\GG}>\tau_{k}$ by the event $T>\tau_{k}$ without
changing the events or probabilities.  Next we can rewrite the second
probability as a sum by conditioning on the variable $L_{k}$, to
obtain that the preceding display is equal to
\begin{eqnarray*}
\lefteqn{\sum_{l_{k}}\Big[
P\left(T>\tau_{k}|\overline{L}_{k-1}=\overline{l}_{k-1},\overline{A}_{k-1}=
\overline{\GG}\left(\overline{l}_{k-1}\right),T>\tau_{k-1}\right)} \\
&&\hspace{1.5cm}
\times P\left(T^{\GG}>t|\overline{L}_{k}=\overline{l}_{k},
\overline{A}_{k-1}=\overline{\GG}\left(\overline{l}_{k-1}\right),
T>\tau_{k}\right)\Big]\\
&&\hspace{1.5cm}
\times P\left(L_{k}=l_{k}|\overline{L}_{k-1}=\overline{l}_{k-1},
\overline{A}_{k-1}=
\overline{\GG}\left(\overline{l}_{k-1}\right),T>\tau_{k}\right).
\end{eqnarray*}
Finally we replace the probability involving the counterfactual
variable $T^\GG$ by the right side of (\ref{Gfc}), which is 
permitted in view of the induction hypothesis.
This yields the right side of (\ref{Gfc}) for $k-1$, 
and concludes the induction step.
\Endproof


\section{Reparameterization}
\label{repars}
To investigate the effect of a given treatment regime $\GG$ on
survival, it suffices to know the conditional distributions given in
(\ref{plm}) and (\ref{ptm}). Given these distributions we can compute
the counterfactual survival functions by using the G--computation
formula, given by Theorem~\ref{Gcomp}.

Because carrying out this computation may be a formidable task, we may
perform the calculation by simulation methods, rather than by
analytical calculation.  Robins~(1986, 1987, 1988) provides a Monte
Carlo algorithm, called the ``Monte Carlo G--computation algorithm'',
for evaluating the functions $s_{\GG}$ that satisfactorily resolves
potential difficulties with the analytical computation. We refer the
reader to these papers for further discussion.

A difficulty is that the distributions in (\ref{plm}) and
(\ref{ptm}) will typically be unknown and must be estimated from the
data. One possibility is to specify models for (\ref{plm}) and
(\ref{ptm}), for instance logistic or Cox models, and next estimate
the unknown parameters from the data. The function $s_\GG$ can then be
estimated using the Monte Carlo G--computation algorithm with model
derived estimates. Robins~(1986, 1987) provides several worked
examples of this approach.

This approach has a number of unattractive features. Estimation of the
function $s_\GG$ according to the preceding scheme and without
confidence intervals, may be feasible, but testing whether treatment
affects the outcome is complicated.  The models used to specify
$s_\GG$ will usually be rough approximations, and the null hypothesis
of no treatment effect will be a complex function of all parameters.
Standard statistical software may not apply, and in large datasets the
null hypothesis will usually be rejected, just because of model
misspecification (cf.\ Robins~(1986, 1987, 1988, 1989)). In this paper
we take a different approach, based on a reparameterization of the
joint distribution of the observations
$\left(\overline{L}_K,\overline{A}_K,T\right)$ using \emph{structural
nested failure time models (SNFTM)}.

SNFTMs are models for the causal effect of skipping a ``last''
treatment dose given the past, thus reverting to the ``baseline
treatment''.  To make this precise, suppose that there is a certain
baseline treatment regime, which we shall refer to as ``no
treatment''. This could for instance be ``zero medication'', and
consequently we shall let a zero in the sets $\overline{\mathcal A}_k$ of
treatment dosages refer to treatment under the baseline treatment
regime.

At any time point $\tau_k$ a doctor could switch a patient to the
baseline regime, at least conceptually, and leave her there. Let
$\left(\overline{a}_k,\overline{0}\right)$ be an abbreviation for the
treatment regime $\GG=\left(a_0,\ldots,a_k,0,\ldots,0\right)$, i.e.\
the $m$th coordinate function of $\GG$ is given by
\begin{equation*}
g_m\left(\overline{l}_m\right)= \left\{\begin{array}{ll} a_m & {\rm
for} \;{\rm any}\;{\rm value}\;{\rm of}\; {\rm the}\; {\rm
covariate}\;{\rm vector}\;\overline{l}_m \; {\rm if}\; m\leq k,\\
0 & {\rm if}\; m>k. \end{array}\right.
\end{equation*}
Henceforth, we shall always assume that Assumptions~\ref{cons}
(consistency) and~\ref{ass} (no unmeasured confounding) are
satisfied. Then, by Theorem~\ref{Gcomp}, if the treatment regime
$(\overline a_k,\overline 0)$ is evaluable, the function 
$$t\mapsto s_{\overline{l}_k,\left(\overline{a}_k,\overline 0\right)}(t)$$
(by
definition the right side of (\ref{Gfc}) with $\GG=(\overline
a_k,\overline 0)$) is the conditional survival function of the
counterfactual survival time
$T^{\left(\overline{a}_k,\overline{0}\right)}$ given the treatment-
and covariate history $\left(\overline{l}_k,\overline{a}_{k-1}\right)$
up to time $\tau_k$, and given that
$T^{\left(\overline{a}_k,\overline{0}\right)}>\tau_k$.  Define
``shift-functions'' $\gamma$ by
\begin{equation}\label{gd}
\gamma_{\overline{l}_k,\overline{a}_k}(t)=
s^{-1}_{\overline{l}_k,\left(\overline{a}_{k-1},\overline{0}\right)}
\circ s_{\overline{l}_k,\left(\overline{a}_{k},\overline{0}\right)} (t),
\end{equation}
where the inverse $s^{-1}$ is the quantile function of the corresponding
survival function.

The functions $\gamma$ map percentiles of the distribution of the
random variable $T^{\left(\overline{a}_k,\overline{0}\right)}$ into
those of the distribution of the random variable
$T^{\left(\overline{a}_{k-1},\overline{0}\right)}$,
\begin{equation}
\label{gdalt}
s_{\overline{l}_k,\left(\overline{a}_{k-1},\overline{0}\right)}
\circ \gamma_{\overline{l}_k,\overline{a}_k}
=s_{\overline{l}_k,\left(\overline{a}_{k},\overline{0}\right)}.
\end{equation}
The functions $\gamma$ thus measure the effect of skipping the ``last''
treatment dose $a_k$ given the covariate and treatment history
$(\overline l_k,\overline a_{k-1})$. We assume that the survival functions
are continuous and strictly decreasing, so that (\ref{gd}) and (\ref{gdalt})
give equivalent definitions.

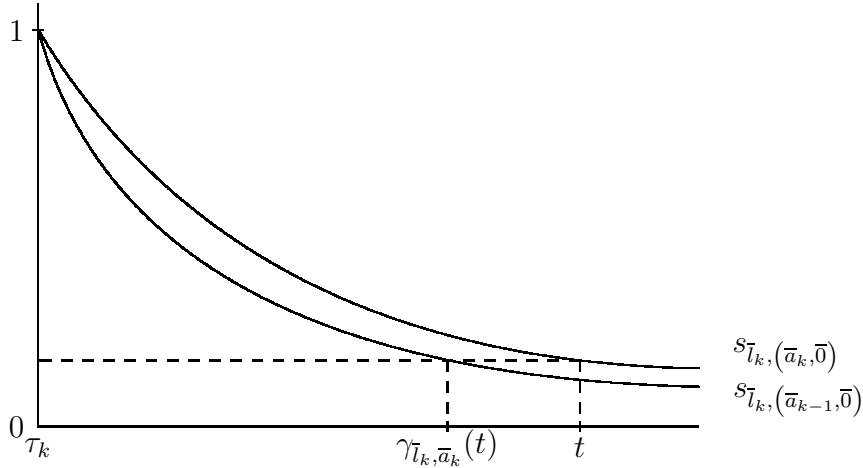
\begin{figure}[htb!]
\begin{picture}(350,170)
\put(40,10){\line(1,0){250}}
\put(40,10){\line(0,1){160}}
\put(32,10){\makebox(0,0){$0$}}
\put(32,160){\makebox(0,0){$1$}}
\put(38,160){\line(1,0){4}}
\put(40,2){\makebox(0,0){$\tau_k$}}
\qbezier[4000](40,160)(75,30)(290,25)
\qbezier[4000](40,160)(115,35)(290,32)
\put(245,2){\makebox(0,0){$t$}}
\put(245,8){\line(0,1){4}}
\put(195,2){\makebox(0,0){$\gamma_{\overline{l}_k,\overline{a}_k}(t)$}}
\put(195,8){\line(0,1){4}}
\put(328,20){\makebox(0,0){$s_{\overline{l}_k,\left(\overline{a}_{k-1},
\overline{0}\right)}$}}
\put(323,37){\makebox(0,0){$s_{\overline{l}_k,\left(\overline{a}_k,
\overline{0}\right)}$}}
\put(245,15){\line(0,1){4}}
\put(245,23){\line(0,1){4}}
\put(245,31){\line(0,1){4}}
\put(195,15){\line(0,1){4}}
\put(195,23){\line(0,1){4}}
\put(195,31){\line(0,1){4}}
\put(40,35){\line(1,0){5}}
\put(49,35){\line(1,0){4}}
\put(57,35){\line(1,0){4}}
\put(65,35){\line(1,0){4}}
\put(73,35){\line(1,0){4}}
\put(81,35){\line(1,0){4}}
\put(89,35){\line(1,0){4}}
\put(97,35){\line(1,0){4}}
\put(105,35){\line(1,0){4}}
\put(113,35){\line(1,0){4}}
\put(121,35){\line(1,0){4}}
\put(129,35){\line(1,0){4}}
\put(137,35){\line(1,0){4}}
\put(145,35){\line(1,0){4}}
\put(153,35){\line(1,0){4}}
\put(161,35){\line(1,0){4}}
\put(169,35){\line(1,0){4}}
\put(177,35){\line(1,0){4}}
\put(185,35){\line(1,0){4}}
\put(193,35){\line(1,0){4}}
\put(201,35){\line(1,0){4}}
\put(209,35){\line(1,0){4}}
\put(217,35){\line(1,0){4}}
\put(225,35){\line(1,0){4}}
\put(233,35){\line(1,0){4}}
\put(241,35){\line(1,0){4}}
\end{picture}
\caption{Illustration of the shift-function $\gamma$.
In this picture the 
function $s_{\overline{l}_k,\left(\overline{a}_{k-1},\overline{0}\right)}$
lies to the left of the function 
$s_{\overline{l}_k,\left(\overline{a}_{k},\overline{0}\right)}$, indicating
that skipping the treatment $a_k$ decreases survival for patients
with covariate and treatment history $(\overline{l}_k,\overline{a}_{k-1})$.
In this case the 
function $\gamma_{\overline{l}_k,\overline{a}_k}$ is below the identity.}
\label{gamf}
\end{figure}

If the ``last treatment'' $a_k$ has no effect, then the functions
$s_{\overline{l}_k,\left(\overline{a}_{k-1},\overline{0}\right)}$ and
$s_{\overline{l}_k,\left(\overline{a}_{k},\overline{0}\right)}$ are
identical, and the function $\gamma_{\overline{l}_k,\overline{a}_k}$
is the identity function. More generally, the function
$\gamma_{\overline{l}_k,\overline{a}_k}$ can be seen to measure the
effect of the treatment $a_k$ given in
$\left[\tau_k,\tau_{k+1}\right)$ on (counterfactual) survival.  This
is illustrated in Figure~\ref{gamf}.

Conversely, if the shift function
$\gamma_{\overline{l}_k,\overline{a}_k}$ is equal to the identity
function, then the distribution of the counterfactual variables
$T^{\left(\overline{a}_k,\overline{0}\right)}$ and
$T^{\left(\overline{a}_{k-1},\overline{0}\right)}$ coincide for
patients with past covariate- and treatment history $\overline{l}_k$
and $\overline{a}_{k-1}$.  This suggests that, if
$\gamma_{\overline{l}_k,\overline{a}_k}$ is the identity function for
all values of $\overline{l}_k$, $\overline{a}_k$ and $k$, then
treatment does not affect the outcome of interest: skipping the last
treatment does not affect the outcome of interest, next skipping the
second-last treatment does not affect the outcome of interest,
etcetera.

For a rigorous proof of this conclusion it is necessary that
sufficiently many treatment regimes are evaluable, because the
functions $s_{\overline l_k,g}$ (defined in terms of the distribution
of the observable data by the right side of (\ref{Gfc})) are equal to
the counterfactual survival distributions only if the treatment regime
$g$ is evaluable.  For instance, the treatment regime
$g=\left(\overline{a}_k,\overline{0}\right)$ need not be evaluable for
all $\overline{a}_k$ and hence the distributions of the counterfactual
variables $T^{\left(\overline{a}_k,\overline{0}\right)}$ and/or
$T^{\left(\overline{a}_{k-1},\overline{0}\right)}$ may not be
identifiable from the observed data.  To overcome this difficulty we
assume that the baseline treatment regime $\overline{0}$ is
``admissible''. A treatment regime is called ``admissible'' if in
\emph{every} situation there is a positive probability for this regime
to be implemented in the next step.  As applied to the baseline regime
$\overline 0$, this property takes the form of the following
assumption.

\begin{assu}\label{base} \emph{(admissible baseline treatment regime).}
For each $k$, each
$\overline{l}_k\in\overline{{\cal L}}_k$ and each
$\overline{a}_{k-1}\in\overline{{\cal A}}_{k-1}$,
\begin{equation*} P\left(\overline{L}_k=\overline{l}_k,
\overline{A}_{k-1}=\overline{a}_{k-1}, T>\tau_k\right)>0 \Rightarrow
P\left(\overline{L}_k=\overline{l}_k,\overline{A}_{k-1}=\overline{a}_{k-1},
A_k=0, T>\tau_k\right)>0.
\end{equation*}
\end{assu}

Under this assumption the shift functions
$\gamma_{\overline{l}_k,\overline{a}_k}$ are identifiable for all
values of $(k,\overline{l}_k, \overline{a}_k)$ with
$P\left(\overline{L}_k=\overline{l}_k,\overline{A}_k=
\overline{a}_k,T>\tau_k\right)>0$, and fully characterize the
potential effect of any treatment regime.  This is the content of the
following theorem, whose proof is deferred to Appendix~\ref{appid}.
(As shown in Lok~(2001, Section 2.12), Assumption~\ref{base} can
be avoided if one allows $\overline{0}$ to be a so-called
admissible baseline course of treatment, which may not only depend on
past covariate- but also on past treatment history. Some
admissible baseline course of treatment, which has a positive
probability of occurring after any observed treatment- and covariate
history, always exists.)

\begin{thm}\label{geq}
Under Assumptions~\ref{ass} (no unmeasured confounding), \ref{cons}
(consistency) and \ref{base} (admissible baseline treatment regime),
the distribution of $T^\GG$ is the same under all evaluable treatment
regimes $\GG$ if and only if the shift-function
$\gamma_{\overline{l}_k,\overline{a}_k}$ is the identity for all $(k,
\overline{l}_k, \overline{a}_k)$ with
$P\left(\overline{L}_k=\overline{l}_k,\overline{A}_k=\overline{a}_k,
T>\tau_k\right)>0$.
\end{thm}

It follows that the functions $\gamma_{\overline{l}_k,\overline{a}_k}$
characterize the null hypothesis of no treatment effect.
Because they also possess an easy interpretation in terms of the
effect of a ``last blip'' of treatment, it is attractive 
to model these functions rather than the set of conditional
distributions in (\ref{plm}) and (\ref{ptm}). 
A \emph{structural nested failure time model} is a parametrized
family of functions used to  model the functions
$\gamma_{\overline{l}_k,\overline{a}_k}$. Each of the model
functions is an increasing function on $[\tau_k,\infty)$ (that
can arise as a quantile-distribution function), with the identity
function referring to the absence of the treatment effect.

\begin{figure}
\qquad\resizebox{4in}{!}{\includegraphics{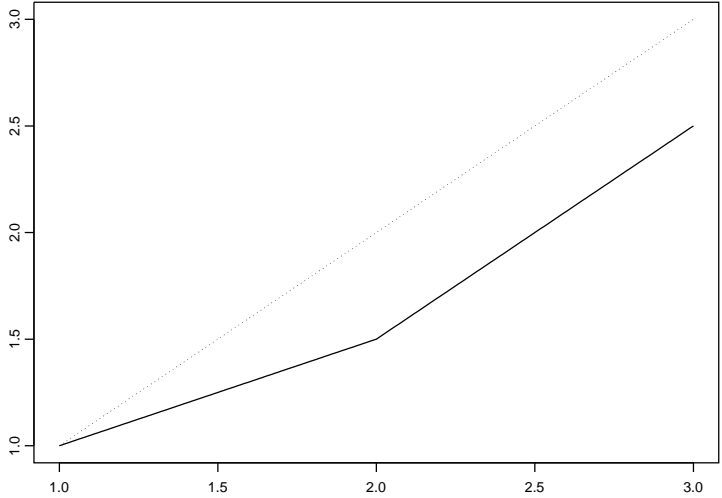}}
\caption{Examples of shift functions. The picture shows the identity
function (dashed) and the function 
$t\mapsto \tau_k
+\left(\min\left\{\tau_{k+1},t\right\}-\tau_k\right)
0.5+\left(t-\tau_{k+1}\right)1_{\left\{t>\tau_{k+1}\right\}}$
for $\tau_k=1<\tau_{k+1}=2$, which corresponds to decreasing
survival by skipping the treatment in the interval $(\tau_k, \tau_{k+1}]$.}
\label{figureSNFTM}
\end{figure}

With the parameter denoted by $\psi=(\psi_1,\psi_2,\psi_3)$,
one example of an SNFTM would be 
\begin{equation*}
\gamma^\psi_{\overline{l}_k,\overline{a}_k}\left(t\right)=\tau_k
+\left(\min\left\{\tau_{k+1},t\right\}-\tau_k\right)
e^{\psi_1 a_k +\psi_2 a_k a_{k-1} +\psi_3 a_k l_k}+
\left(t-\tau_{k+1}\right)1_{\left\{t>\tau_{k+1}\right\}}.
\end{equation*}
If $\psi=0$, then this function reduces to the identity function,
indicating that the parameter value $\psi=0$ corresponds to
the absence of a treatment effect. For nonzero values of $\psi$
the model corresponds to a 
``change of time scale''  depending on present and
past treatment $\left(a_k,a_{k-1}\right)$ and present covariate
($l_k$). The variable $L_k$ might for instance be
the univariate covariate CD4 lymphocyte count at $\tau_k$, and
the variable $A_k$ the AZT prescription. Then the given model
allows for interaction between CD4 lymphocyte count and
treatment, and could of course be extended with other factors.
Figure~\ref{figureSNFTM} shows two typical functions
$\gamma$ following this model.

\section{Mimicking counterfactual outcomes}
\label{mimsec}
In the next two sections we present two methods for estimating the
parameter $\psi$ in a structural nested failure time
model. Theorem~\ref{Tog} below is basic for both methods. Consider the
following transformation of the observation
$\left(\overline{L}_K,\overline{A}_K,T\right)$, using the ``true'' shift
functions $\gamma$ (given by (\ref{gd})):
\begin{equation} \label{Togd}
T_0^\gamma=\gamma_{\overline{L}_0,\overline{A}_0}
\circ \gamma_{\overline{L}_1,\overline{A}_1}
\circ \cdots
\circ \gamma_{\overline{L}_{p\left(T\right)},
\overline{A}_{p\left(T\right)}}\left(T\right),
\end{equation}
where $p(T)=\max\left\{k:\tau_k<T\right\}$.
The application of the function
$\gamma_{\overline{L}_{p\left(T\right)},\overline{A}_{p\left(T\right)}}$
to $T$ annihilates the effect of the last treatment
$A_{p\left(T\right)}$, and each further application of a shift
function annihilates the effect of an earlier treatment.  This
explains the following theorem, which is proved in
Appendix~\ref{appTog}.

\begin{thm} \label{Tog}\emph{(mimicking counterfactual outcomes).}
The variable $T_0^\gamma$ defined in (\ref{Togd})
possesses survival function $s_{\overline{0}}$.
Furthermore, for every $k\geq0$,
\begin{equation}\label{indep}
A_k \cip T_0^\gamma |\overline{L}_k,\overline{A}_{k-1},T>\tau_k.
\end{equation}
\end{thm}

The variable $T_0^\gamma$ is a (deterministic) function of the data
vector $\left(\overline{L}_K,\overline{A}_K,T\right)$, through the
(unknown) family of shift-functions $\gamma$. If the shift functions
$\gamma$ would be known, then we would be able to ``mimic'' the
survival time without treatment by calculating the transformation
$T_0^\gamma$. By the preceding theorem this variable is distributed
according to $s_{\overline 0}$ and hence under the conditions of 
Theorem~\ref{Gcomp} possesses the same distribution as $T^{\GG}$ for
$\GG=\overline{0}$, the null treatment.

Equation (\ref{indep}) shows that the variable $T_0^\gamma$ also
shares the ``no unmeasured confounding'' property
(Assumption~\ref{ass}) of counterfactual variables (in a slightly
stronger form).

\section{Maximum likelihood estimation}
\label{mlesec}
In this section we consider likelihood based inference for the
parameter $\psi$ in a given SNFTM. Clearly this requires that we make
the parameter $\psi$ visible in the density of the observation
$\left(\overline{L}_K,\overline{A}_K,T\right)$. We first show that this
can be achieved using the transformation $T_0^\gamma = T_0^\gamma
\left(T,\overline{L}_K,\overline{A}_K\right)$ defined in
(\ref{Togd}), which will depend on $\psi$ if we use a
SNFTM for $\gamma$.

\begin{thm}\label{mleth}\emph{(the likelihood rewritten).}
Suppose that Assumption~\ref{base} (admissible baseline treatment
regime) holds.
Suppose moreover that $\left(T,\overline{L}_K,\overline{A}_K\right)$ has a
Lebesgue density, and that the function $t\mapsto
s_{\left(\overline{l}_k,\left(\overline{a}_k,\overline{0}\right)\right)}
\left(t\right)$ is continuously differentiable in $t$, for all
$\overline{l}_k$, $\overline{a}_k$ with
$P\left(\overline{L}_k=\overline{l}_k,\overline{A}_k=\overline{a}_k,
T>\tau_k \right)>0$, with strictly negative derivative except for at
most finitely many points. Then the joint density of
$\left(T,\overline{L},\overline{A}\right)$ can be rewritten as
\begin{eqnarray*}
\lefteqn{f_{T,\overline{L},\overline{A}}\left(t,\overline{l},
\overline{a}\right)}\\
&=&\frac{\partial}{\partial t} t_0^\gamma
\left(t,\overline{l}_p,\overline{a}_p\right)
f_{T_0^\gamma}\left(t_0^\gamma\left(t,\overline{l}_p,\overline{a}_p\right)
\right)
P\left(L_0=l_0|T_0^\gamma=t_0^\gamma\right)
P\left(A_0=a_0|L_0=l_0\right)\nonumber\\
&&\prod_{k=0}^p \Big\{P\left(L_k=l_k|\overline{L}_{k-1}=\overline{l}_{k-1},
\overline{A}_{k-1}=\overline{a}_{k-1},T>\tau_k, T_0^\gamma=t_0^\gamma\right)
\\&&\hspace{3cm}
P\left(A_k=a_k|\overline{L}_{k}=\overline{l}_{k},
\overline{A}_{k-1}=\overline{a}_{k-1},T>\tau_k\right)\Big\},
\end{eqnarray*}
where $\tau_p<t\le \tau_{p+1}$ and
\begin{equation*}
t_0^\gamma\left(t,\overline{l}_p,\overline{a}_p\right)=
\gamma_{\overline{l}_0,\overline{a}_0}
\circ \gamma_{\overline{l}_1,\overline{a}_1}
\circ \cdots
\circ \gamma_{\overline{l}_{p},
\overline{a}_{p}}\left(t\right).
\end{equation*}
\end{thm}

\proof
Under the conditions of Theorem~\ref{mleth},
\begin{equation*}
\left(T,\overline{L},\overline{A}\right)\mapsto
\left(T_0^\gamma,\overline{L},\overline{A}\right)=
\left(t_0^\gamma\left(T\right),\overline{L},\overline{A}\right)
\end{equation*}
is a one-to-one mapping. Thus if $t_0^\gamma$ were continuously
differentiable everywhere, then the identity
\begin{equation}\label{lik1}
f_{T,\overline{L},\overline{A}}\left(t,\overline{l},\overline{a}\right)
=\frac{\partial}{\partial t} t_0^\gamma
\left(t,\overline{l},\overline{a}\right)
f_{T_0^\gamma,\overline{L},\overline{A}}
\left(t_0^\gamma\left(t,\overline{l},\overline{a}\right),\overline{l},
\overline{a}\right)
\end{equation}
would be immediate from the change of variables formula.  We show that
(\ref{lik1}) holds under the conditions of Theorem~\ref{mleth}
too. Next the assertion of the theorem follows by repeated
conditioning and using Theorem~\ref{Tog}.

To prove (\ref{lik1}) in general, note that the probability space
consists of countably many sets of the form
$\left(\overline{L}_K=\overline{l}_K,\overline{A}_K=\overline{a}_K\right)$, so
that by countable additivity of measures it suffices to prove
(\ref{lik1}) on each of these sets that has probability greater than
$0$. On each of these sets, $t_0^\gamma$ is one-to-one and
continuously differentiable except for at finitely many points: it is
the composition of finitely many functions
$\gamma_{\overline{l}_k,\overline{a}_k}$ and under the assumptions of
Theorem~\ref{mleth},
\begin{equation*}\gamma'_{\overline{l}_k,\overline{a}_k}\left(t\right)
=\bigl(s^{-1}_{\overline{l}_k,\left(\overline{a}_{k-1},\overline{0}\right)}
\circ 
s_{\overline{l}_k,\left(\overline{a}_{k},\overline{0}\right)}\bigr)'(t)
=\frac{1}{s'_{\overline{l}_k,\left(\overline{a}_{k-1},\overline{0}\right)}
\bigl(\gamma_{\overline{l}_k,\overline{a}_k}\left(t\right)\bigr)}
s'_{\overline{l}_k,\left(\overline{a}_{k},\overline{0}\right)}\left(t\right)
\end{equation*}
exists and is continuous except for at most finitely many $t$.
Thus, from the change of variables formula, equation (\ref{lik1}) is
true on each set
$\left(\overline{L}_K=\overline{l}_K,\overline{A}_K=\overline{a}_K\right)$, as
we needed to show. 
\Endproof

Regarding the conditions of Theorem~\ref{mleth} we note that the
baseline treatment regime $\overline{0}$ may not be constant, whence
the death rate under $\overline{0}$ may change at the time points
$\tau_m$. However, it will often be reasonable to assume
differentiability of the function
$s_{\left(\overline{l}_k,\left(\overline{a}_k,\overline{0}\right)\right)}
\left(t\right)$ on all intervals $\left(\tau_m,\tau_{m+1}\right)$.

For likelihood inference concerning the parameter $\psi$ of an SNFTM, we
shall generally drop the factors
\begin{equation}
\label{condlawtreatment}
P\left(A_k=a_k|\overline{L}_k=\overline{l}_k,
\overline{A}_{k-1}=\overline{a}_{k-1}\right)
\end{equation}
from the likelihood. All other terms involve $\psi$ through
$T_0^\gamma$ and we will need to specify models for these terms in
order to proceed, typically involving additional parameters.  Given
such models we can estimate $\psi$ by the corresponding coordinate of
the maximum likelihood estimator obtained by maximizing the likelihood
over all parameters. Of course finding this maximizer may be a
formidable task.

Since the null hypothesis of no treatment effect is equivalent to the
functions $\gamma_{\overline{l}_{k},\overline{a}_k}$ being equal to
the identity function, by Theorem~\ref{geq}, this hypothesis can be
fully expressed in the parameter $\psi$. For instance, we could, by
convention, construct our SNFTM in such a way that this null
hypothesis is equivalent to $H_0: \psi=0$. Then we can obtain a
likelihood-based test for the null hypothesis of no treatment effect
using the Wald, score or likelihood ratio test for $H_0: \psi=0$.

\section{G--estimation}
\label{Gest}
The likelihood methods of the preceding section require the 
specification of models for the conditional laws of the covariates,
among others, next to a specification of an SNFTM. In this section we
present an alternative approach to testing and estimation of the
parameter in a SNFTM, called \emph{G--estimation} in
Robins~(1998). This approach is based on models for the conditional
distributions of the treatment variables given in
(\ref{condlawtreatment}). It can be considered a semiparametric
approach, where the parametric component refers to the laws
(\ref{condlawtreatment}) and all other laws appearing in
Theorem~\ref{mleth} form the nonparametric, unspecified
component. From a practical perspective modelling the distributions
(\ref{condlawtreatment}) is more attractive than modelling the
remaining laws in Theorem~\ref{mleth}, as it may be expected that
doctors have clear ideas, at least qualitatively, about how they reach
their decisions about treatment.

The method of G--estimation is based on the conditional independence
of the ``blipped-up'' variable $T_0^\gamma$ defined in (\ref{Togd})
and the treatment variable $A_k$ given the variables $\overline L_k$
and $\overline A_{k-1}$, for each $k$, asserted by
Theorem~\ref{Tog}. Consider first testing the null hypothesis $H_0:
\gamma=\gamma_0$ for a given shift function $\gamma_0$. 
Theorem~\ref{Tog} gives, under the null hypothesis, that, for each $k$,
\begin{equation}\label{ind}
A_k \cip T_0^{\gamma_0}\,|\,\overline{L}_k,\overline{A}_{k-1},T>\tau_k.
\end{equation}
This is an assertion about the observed data vector
$(\overline L_K,\overline A_K, T)$ only. Any test for the validity
of (\ref{ind}) is therefore a test for the null hypothesis 
$H_0: \gamma=\gamma_0$.

In order to operationalize this idea we adopt for each $k$ a model
\begin{equation*}
P_\theta\left(A_k=a_k|\overline{L}_k=\overline{l}_k,
\overline{A}_{k-1}=\overline{a}_{k-1},T>\tau_k\right)
\end{equation*}
for the prediction of treatment given the past, indexed by some
parameter $\theta$. Such a model tries to explain the treatment $A_k$
by the values of the covariates up to time $\tau_k$ and the preceding
treatment history.  Formula (\ref{ind}) implies that, under the null
hypothesis, inclusion of the variable $T_0^{\gamma_0}$ as an extra
explanatory variable is useless for the prediction of $A_k$, if past
covariate- and treatment information $\overline{L}_k$ and
$\overline{A}_{k-1}$ are known. Thus given a term of the form
$\alpha\, T_0^{\gamma_0}$ in the prediction model with $\alpha$ a
parameter, the true value of $\alpha$ must be equal to $0$, because of
(\ref{ind}). It follows that we can test the null hypothesis $H_0:
\gamma=\gamma_0$ by adding a term $\alpha T_0^{\gamma_0}$ anyway, and
next test the null hypothesis $H_0: \alpha=0$ in the model indexed by
the overall parameter $(\theta,\alpha)$.  Depending on the chosen
types of model such a test, for instance a Wald, score or the
likelihood ratio test, can be performed by standard statistical
software.

This procedure is particularly simple for testing the
null hypothesis of no treatment effect. In view of 
Theorem~\ref{geq}, this is equivalent to testing whether the function
$\gamma$ is equal to the identity function, i.e.\ we take $\gamma_0$ in the
preceding equal to the identity function. In this case
the variable $T_0^{\gamma_0}$ is equal to $T$, and hence
the G--estimation procedure reduces to testing the null hypothesis
$H_0: \alpha=0$ in a regression model that tries to explain the variable
$A_k$ by the variables $\overline L_k$, $\overline A_{k-1}$ and $\alpha T$. 
The null hypothesis of no treatment effect can be tested
in this way without specifying a model for the shift function $\gamma$.

For a specific example, suppose that the treatment variables $A_k$ are
binary-valued.  Then a logistic regression model is a standard choice
for modelling the probabilities (\ref{condlawtreatment}).  We might
add the variable $\alpha T_0^{\gamma}$ to a logistic regression model
to form the model
\begin{equation*}
P_{\theta,\alpha}\left(A_k=a_k|\overline{L}_k,\overline{A}_{k-1},T>\tau_k,T_0^{\gamma}\right)
=\frac{1}{1+e^{\theta\cdot f_k(\overline{L}_k,\overline{A}_{k-1})
+\alpha g_k(T_0^{\gamma})}},
\end{equation*}
for given, known functions $f_k$ and $g_k$, and unknown parameters
$\theta$ and $\alpha$. A test for the null hypothesis $H_0: \alpha=0$
can be carried out by standard software for logistic regression.

Given an SNFTM $\psi\mapsto \gamma_\psi$ for the shift functions
$\gamma$, indexed by a parameter $\psi$, we can extend the preceding
testing methods to full inference on the parameter $\psi$.  First, we
can obtain confidence regions for $\psi$ by inverting the tests for
the null hypotheses $H_0: \gamma=\gamma_{\psi}$ in the usual way: the
value $\psi$ belongs to the confidence region if the corresponding
null hypothesis $H_0$ is not rejected.

A natural estimator of $\psi$ would be the center of a confidence
set, or, alternatively, a value of $\psi$ for
which $T_0^{\gamma_\psi}$ contributes the
least to the prediction model for treatment given the
past. That is, the $\psi$ for which the fitted model for
\begin{equation}\label{akpred}
P_{\theta,\alpha}\left(A_k=a_k|\overline{L}_k,
\overline{A}_{k-1},T>\tau_k,\alpha T_0^{\gamma_\psi}\right).
\end{equation}
does not include the variable $T_0^{\gamma_\psi}$, i.e.\ $\alpha=0$.
For each given value of the parameter $\psi$ of the SNFTM we may
obtain estimators $\hat{\theta}(\psi)$ and $\hat{\alpha}(\psi)$ for
the parameters $\theta$ and $\alpha$, based on the observations
$(\overline{L}_K^i,\overline{A}_K^i,T^i)$ on $n$ persons.  Then we
define $\hat{\psi}$ as the solution of the equation
\begin{equation*}
\hat{\alpha}\left(\psi\right)=0.
\end{equation*}
If we use a logistic regression model, then the estimators
$\hat\theta$  and $\hat\alpha$ can be obtained with
standard software, for each given value of $\psi$.
The estimator $\hat\psi$ can next be found by a grid search method.
Alternatively, we can implement a direct numerical method for
estimating $\psi$.

The procedures just outlined may appear a bit unusual, in view of
their indirect nature. However, in most cases they can also be
interpreted in a standard way. For instance, the procedure for
estimating $\alpha$ for given $\psi$ will often be equivalent to
solving $\hat{\alpha}=\hat{\alpha}\left(\psi\right)$ from an
estimating equation of the type
\begin{equation*}
\sum_{i=1}^n h_{\alpha,\psi}\bigl(\overline{L}_K^i,\overline{A}_K^i,T^i\bigr)
=0.
\end{equation*}
Then $\hat{\psi}$ satisfying
$\hat{\alpha}\bigl(\hat{\psi}\bigr)=0$ will satisfy
the estimating equation
\begin{equation*}
\sum_{i=1}^n h_{0,\hat{\psi}}\bigl(\overline{L}_K^i,\overline{A}_K^i,T^i\bigr)
=0.
\end{equation*}
Because $\alpha(\psi_0)=0$ for the true value $\psi_0$ of $\psi$,
the true value of $\psi$ is a solution to the equation
\begin{equation*}
E h_{0,\psi}\left(\overline{L}_K,\overline{A}_K,T\right)=0.
\end{equation*}
In other words, $\hat{\psi}$ will be the solution of an unbiased
estimating equation, whence the (asymptotic) properties of $\hat\psi$
can be ascertained with the usual theory for M-estimators
(e.g.\ Van der Vaart~(1998)). For
instance, we may expect the sequence
$\sqrt{n}\bigl(\hat{\psi}-\psi\bigr)$ to be asymptotically (as
$n\rightarrow\infty$) normal with mean zero and variance
\begin{equation*}
\frac{E h_{0,\psi}^2\left(\overline{L}_K,\overline{A}_K,T\right)}
{\bigl(\frac{\partial}{\partial \psi} E 
h_{0,\psi}\left(\overline{L}_K,\overline{A}_K,T\right)\bigr)^2}.
\end{equation*}
Lok~(1991) has studied the validity of these results, and has thus
justified the preceding procedures.

\section{Summary and extensions}
We have shown that the AZT treatment regime-specific, counterfactual
AIDS-free survival curves $P\left(T^g>t\right)$ are identified for all
evaluable treatment regimes $g$ if our maintained assumption of no
unmeasured confounding, Assumption~\ref{ass}, is met. This assumption
will hold if the investigator has succeeded in recording in
$\overline{l}_k$ data on all covariates that, conditional on past AZT
history $\overline{a}_{k-1}$, predict both the AZT dosage rate $a_k$
in $\left(\tau_k,\tau_{k+1}\right]$ and the random variables $T^g$
representing time to AIDS had, contrary to fact, all subjects followed
an AZT treatment history consistent with regime $g$.

Further, we have shown that, under the assumption of no unmeasured
confounding, Assumption~\ref{ass}, the shift functions $\gamma$ of an
SNFTM are the identity function if and only if the G--null hypothesis
of no causal effect of AZT on time to AIDS is true. We have expressed
the likelihood of the observable random variables
$\left(T,\overline{L}_K,\overline{A}_K\right)$ in terms of the transformed
random variables $\left(T_0^\gamma,
\overline{L}_K,\overline{A}_K\right)$. We then developed parametric
likelihood based tests of the hypothesis $\gamma={\rm id}$ by specifying fully
parametric models for the joint distribution of
$\left(T,\overline{L}_K,\overline{A}_K\right)$ in terms of the
transformed random variables $\left(T_0^\gamma,
\overline{L}_K,\overline{A}_K\right)$.

Even in the absence of censoring or missing data, a major limitation
of the fully parametric likelihood-based tests of the null hypothesis
$\gamma={\rm id}$ from Section~\ref{mlesec} is that misspecification of the
parametric models for the distribution of $L_k$ given
$\overline{L}_{k-1}$, $\overline{A}_{k-1}$ and $T_0^\gamma$, or for
the distribution of $T^{\overline{0}}$, can cause the true
$\alpha$-level of the test to deviate from the nominal
$\alpha$-level. This limitation raised the question of whether it is
possible to construct $\alpha$-level tests of the null hypothesis
$\gamma={\rm id}$ and of more general hypotheses concerning $\gamma$, which
are asymptotically distribution-free. A closely related question is
whether there exist $n^{1/2}$-consistent asymptotically normal
estimators of the parameter $\psi$ of a correctly specified structural
nested failure time model if the joint distribution of the observables
$\left(\overline{L}_K,\overline{A}_K,T\right)$ is otherwise unspecified,
i.e.\ if the distribution of $L_k$ given $\overline{L}_{k-1}$,
$\overline{A}_{k-1}$ and $T_0^\gamma$ and the distribution of the variable
$T^{\overline{0}}$ are left completely unspecified. In
Section~\ref{Gest} we showed that one only needs to specify a parametric
model for the shift function $\gamma$, which models the causal effect
of one treatment dosage given the past, and a parametric model for the
distribution of actual treatment dosage given past treatment- and
covariate history. Doctors will usually have clear ideas about this
latter distribution of treatment decisions. Moreover, the doctors'
interest will often be in the causal effect of one treatment dosage
given the past.

If the null hypothesis of no treatment effect has been rejected and
the parameter $\psi$ of the shift function $\gamma$ has been
estimated, one might wish to estimate the survival distribution
$t\mapsto P\bigl(T^\GG>t\bigr)$ of the outcome under specific
treatment regimes $\GG$ in a way consistent with the estimator
$\hat{\psi}$. This can be done by estimating the distribution of
$T^{\overline{0}}$ (e.g.\ by the empirical distribution of
$T_0^{\gamma^\psi}$) and the empirical distribution of $L_k$ given
$\overline{L}_{k-1}$, $\overline{A}_{k-1}$ and $T_0^\gamma$
($k=0,\ldots,K$) for histories $\overline{L}_{k-1}$,
$\overline{A}_{k-1}$ consistent with $\GG$.  An approximate sample
$\tilde{T}^\GG_i$ ($i=1,2,\ldots$) from the distribution of $T^\GG$
could then be generated by using these estimated distributions: first
draw $T'_0$ from the distribution of $T^{\overline{0}}$, then draw
$L'_0$ from the distribution of $L_0$ given $T_0^\gamma=T'_0$, then
put $A'_0=\GG\bigl(L'_0\bigr)$, then draw $L'_1$ from the distribution
of $L_1$ given $T_0^\gamma=T'_0$, $A_0=A'_0$ and $L_0=L'_0$, etcetera.
Finally put
\begin{equation*}\tilde{T}^\GG=
{\gamma^{\hat{\psi}}_{\overline{L}'_K,\overline{A}'_K}}^{-1}\circ\ldots\circ
{\gamma^{\hat{\psi}}_{\overline{L}'_0,\overline{A}'_0}}^{-1}\bigl(T'_0\bigr).
\end{equation*}
This variable will be generated from the desired distribution.

Extensions of the results of this paper that allow for censoring and
missing data are discussed in Robins~(1988, 1992, 1993, 1998), and Robins
et al~(1992). The extension of G--tests and estimators to continuous
$L_k$ and $A_k$ are discussed in Robins~(1992, 1993), Robins et
al.~(1992), and Gill and Robins~(2001). Robins~(1998) and Lok~(2001)
show that the results in this paper can be extended to allow for jumps
in the treatment- and covariate processes in continuous time.

\begin{appendix}

\section{Alternative formulation of the null hypothesis}
\label{appid}
In this appendix we prove Theorem~\ref{geq} through two lemmas.  The
first lemma shows that if all functions $\gamma$ are equal to the
identity function, then all survival curves $P\left(T^g>t\right)$ for
evaluable treatment regimes are the same. The second lemma shows the
reverse.

\begin{lem}
\label{alem1}
Suppose that Assumptions~\ref{ass} (no unmeasured confounding),
\ref{cons} (consistency) and~\ref{base} (admissible baseline treatment
regime) hold. If $\gamma_{\overline{l}_k,\overline{a}_k}$ is the
identity function for all $k$, $\overline{l}_k\in\overline{{\cal
L}}_k$ and $\overline{a}_k\in\overline{A}_k$ with
$P\left(\overline{L}_k=\overline{l}_k,\overline{A}_k=\overline{a}_k,
T>\tau_k\right)>0$, then
all survival curves $P\left(T^\GG>t\right)$ for evaluable treatment
regimes $\GG$ are the same.
\end{lem}

\proof
We show that for all evaluable treatment regimes $\GG$ and
all $\overline{l}_k$ with \linebreak
$P\left(\overline{L}_k=\overline{l}_k,\overline{A}_k=
\overline{\GG}\left(\overline{l}_k\right),T>\tau_k\right)>0$,
the conditional distributions of the counterfactual variables
$T^\GG$ and 
$T^{\left(\overline{\GG}_{k-1}\left(\overline{l}_{k-1}\right),
\overline{0}\right)}$
given $\overline{L}_k=\overline{l}_k,\overline{A}_{k-1}
=\overline{\GG}\left(\overline{l}_{k-1}\right),T>\tau_k$
are the same, i.e., for $t\ge \tau_k$,
\begin{equation}\label{kweg}
s_{\overline{l}_k,\GG}(t)=
s_{\overline{l}_k,\left(\overline{\GG}_{k-1}\left(\overline{l}_{k-1}\right),
\overline{0}\right)}(t).
\end{equation}
For $k=-1$ this should be read as $s_\GG\left(t\right)=s_{\overline{0}}(t)$,
which implies Lemma~\ref{alem1}.

We prove (\ref{kweg}) by backward induction on $k$, for $t$ fixed.
With $\tau_p$ the last clinic visit time strictly before $t$, we start
with $k=p$ and end with $k=0$. The statement for $k=-1$ follows from
the statement for $k=0$ by summation over $l_0$.

Basis: For $k=p$, by the definition of $s$ as the right side
of (\ref{Gfc}),
$$s_{\overline{l}_p,\GG}(t)
=P\left(T>t|\overline{L}_p=\overline{l}_p,\overline{A}_p=
\overline{\GG}_p\left(\overline{l}_p\right),T>\tau_p\right)
=s_{\overline{l}_p,\left(\overline{\GG}_p\left(\overline{l}_p\right),
\overline{0}\right)}(t),$$
by another application of the definition of $s$.
The right side is equal to 
$s_{\overline{l}_p,\left(\overline{\GG}_{p-1}\left(\overline{l}_{p-1}\right),
\overline{0}\right)}(t)$ by the assumption that
the function $\gamma_{\overline{l}_p,\overline{a}_p}$ with 
$\overline a_p=\overline g_p(\overline l_p)$, is the identity function
is the identity.

Induction step: we suppose that (\ref{kweg}) is true for $k\ge 1$ and
establish (\ref{kweg}) for $k-1$. By straightforward algebra using the
definition of $s_{\overline l_{k-1},g}$,
\begin{eqnarray*} s_{\overline{l}_{k-1},\GG}\left(t\right)
&=& P\left(T>\tau_{k}|\overline{L}_{k-1}=\overline{l}_{k-1},
\overline{A}_{k-1}=\overline{\GG}\left(\overline{l}_{k-1}\right),
T>\tau_{k-1}\right)\\
&&\hspace{0.5cm}\sum_{l_{k}}
P\left(L_{k}=l_{k}|\overline{L}_{k-1}=\overline{l}_{k-1},
\overline{A}_{k-1}=\overline{\GG}\left(\overline{l}_{k-1}\right),
T>\tau_{k}\right) s_{\overline{l}_{k},\GG}\left(t\right).
\end{eqnarray*}
Here we can replace $s_{\overline l_{k},g}$ using the induction
hypothesis, giving that the preceding display is equal to
\begin{eqnarray*}
&&P\left(T>\tau_{k}|\overline{L}_{k-1}=\overline{l}_{k-1},
\overline{A}_{k-1}=\overline{\GG}\left(\overline{l}_{k-1}\right),
T>\tau_{k-1}\right)\\
&&\hspace{0.5cm}
\sum_{l_{k}}
P\left(L_{k}=l_{k}|\overline{L}_{k-1}=\overline{l}_{k-1},
\overline{A}_{k-1}=\overline{\GG}\left(\overline{l}_{k-1}\right),
T>\tau_{k}\right)
s_{\overline{l}_{k},\left(\overline{\GG}_{k-1}\left(\overline{l}_{k-1}\right),
\overline{0}\right)}\left(t\right)\\
&&\hspace{1cm}= s_{\overline{l}_{k-1},\left(\overline{\GG}_{k-1}\left(\overline{l}_{k-1}
\right), \overline{0}\right)}\left(t\right)\\
&&\hspace{1cm}= s_{\overline{l}_{k-1},\left(\overline{\GG}_{k-2}\left(\overline{l}_{k-2}
\right),\overline{0}\right)}\left(t\right),
\end{eqnarray*}
where we use the definition of $s$ in the first equality, and the
assumption that $\gamma_{\overline{l}_{k-1},\overline{a}_{k-1}}$,
for  $\overline{a}_{k-1}=\overline{\GG}_{k-1}(\overline{l}_{k-1})$, is the
identity function in the second.  
\Endproof

\begin{lem}\label{alem2}
Suppose that Assumptions~\ref{ass} (no unmeasured confounding),
\ref{cons} (consistency) and~\ref{base} (admissible baseline treatment
regime) hold. If the survival curves
$P\left(T^\GG>t\right)$ are the same for all evaluable treatment regimes
$\GG$, then the shift function
$\gamma_{\overline{l}_k,\overline{a}_k}$ is the identity for all
$k$, $\overline{l}_k\in\overline{{\cal L}}_k$ and
$\overline{a}_k\in\overline{A}_k$ with
$P\left(\overline{L}_k=\overline{l}_k,\overline{A}_k=\overline{a}_k,
T>\tau_k\right)>0$.
\end{lem}

\proof
Let fixed $\overline{l}_k$, $\overline{a}_k$ with
$P\left(\overline{L}_k=\overline{l}_k,\overline{A}_k=\overline{a}_k,
T>\tau_k\right)>0$ be given. To prove that
$\gamma_{\overline{l}_k,\overline{a}_k}$ is the identity we need to
show that, for all $t>\tau_k$,
\begin{equation}\label{akweg}
s_{\overline{l}_k,\left(\overline{a}_k,\overline{0}\right)}(t)=
s_{\overline{l}_k,\left(\overline{a}_{k-1},\overline{0}\right)}(t).
\end{equation}
Define a treatment regime $\GG^1$ by the coordinate functions
$\GG_m^1\bigl(\overline{\tilde{l}}_m\bigr)=a_m$ if
$\overline{\tilde{l}}_m$ is the initial part of $\overline{l}_k$,
and by $\GG_m^1\bigl(\overline{\tilde{l}}_m\bigr)=0$ otherwise. Define
a second treatment regime $\GG^2$ by 
and $\GG^2=\bigl(\overline{\GG^1}_{k-1},\overline{0}\bigr)$.
Because of Assumption~\ref{base} and because 
$P\left(\overline{L}_k=\overline{l}_k,\overline{A}_k=\overline{a}_k,
T>\tau_k\right)>0$, the treatment
regimes $\GG^1$ and $\GG^2$ are evaluable.   Thus, by assumption,
we have that $P\left(T^{\GG_1}>t\right)=P\left(T^{\GG_2}>t\right)$, and these
probabilities are given by the G--computation formula,
given in Theorem~\ref{Gcomp}. For the first regime this formula
can be written in the form
\begin{eqnarray*}
\lefteqn{P\left(T^{\GG_1}>t\right)}\\
&=&
\sum_{\tilde{l}_0}\cdots
\sum_{\tilde{l}_k} 1_{\overline{\tilde{l}}_k\neq\overline{l}_k}
\prod_{m=0}^k\Big\{ P\bigl(T>\tau_m|\overline{L}_{m-1}=
\overline{\tilde{l}}_{m-1},
\overline{A}_{m-1}=\overline{\GG^1}\bigl(\overline{\tilde{l}}_{m-1}\bigr),
T>\tau_{m-1}\bigr)\\
&&\hspace{1.5cm}
P\bigl(L_m=\tilde{l}_m| \overline{L}_{m-1}=\overline{\tilde{l}}_{m-1},
\overline{A}_{m-1}=\overline{\GG^1}\bigl(\overline{\tilde{l}}_{m-1}\bigr),
T>\tau_{m}\bigr)
\Big\}s_{\overline{\tilde{l}}_k,\GG^1}(t)\\
&&+\bigg[\prod_{m=0}^k
\Big\{P\bigl(T>\tau_m|\overline{L}_{m-1}=\overline{l}_{m-1},
\overline{A}_{m-1}=\overline{\GG^1}\left(\overline{l}_{m-1}\right),
T>\tau_{m-1}\bigr)\\
&&\hspace{1.5cm}
P\bigl(L_m=l_m| \overline{L}_{m-1}=\overline{l}_{m-1},
\overline{A}_{m-1}=\overline{\GG^1}\bigl(\overline{l}_{m-1}\bigr),
T>\tau_{m}\bigr)\Big\}\bigg] s_{\overline{l}_k,\GG^1}(t).
\end{eqnarray*}
A similar expression holds for the treatment regime $\GG^2$.  Because
the regimes $\GG^1$ and $\GG^2$ are constructed to be the same up to
time $\tau_{k-1}$, only the second terms of the summs differs between these
two expressions. Even there, the product preceding
$s_{\overline{l}_k,\GG^1}(t)$ and $s_{\overline{l}_k,\GG^2}(t)$ is the
same for $\GG^1$ and $\GG^2$. Moreover, this factor is strictly
positive, since
$P\left(\overline{L}_k=\overline{l}_k,\overline{A}_k=\overline{a}_k,
T>\tau_k\right)>0$ by assumption. The equality of
$P\left(T^{\GG_1}>t\right)$ and $P\left(T^{\GG_2}>t\right)$ therefore
implies the equality of $s_{\overline{l}_k,\GG^1}(t)$ and
$s_{\overline{l}_k,\GG^2}(t)$. By construction of $\GG^1$ and $\GG^2$,
equation (\ref{akweg}) and hence Lemma~\ref{alem2} follow.  
\Endproof

\section{Mimicking counterfactual outcomes} \label{appTog}
For $t>0$ define $p(t)$ by $\tau_{p(t)}<t\le \tau_{p(t)+1}$,
i.e.\ $\tau_{p(t)}$ is the last clinic visit time strictly before $t$.
For $k\geq 0$ with $k\le p(T)$ we define a random variable by 
\begin{equation*}
T_k^\gamma=\gamma_{\overline{L}_k,\overline{A}_k}\circ
\cdots\circ\gamma_{\overline{L}_p(T),\overline{A}_{p(T)}}(T).
\end{equation*}
For $k>p(T)$ we interprete the (empty) composition of transformations
on the right as the identity and define $T_k^\gamma=T$.

In this appendix we prove the following theorem, which generalizes the
first part of Theorem~\ref{Tog}. This theorem implies the second part,
since $T_0^\gamma$ is a function of
$\bigl(\overline{L}_{k-1},\overline{A}_{k-1},T_k^\gamma\bigr)$.

\begin{thm}
For $t>\tau_k$ and every $\overline l_k$, $\overline a_k$ with
$P\left(\overline{L}_k=\overline{l}_k,\overline{A}_k=
\overline{a}_k,T>\tau_k\right)>0$,
\begin{eqnarray*}
P\left(T^\gamma_k>t|\overline{L}_k=\overline{l}_k,\overline{A}_k=
\overline{a}_k,T>\tau_k\right)
&=&P\left(T^\gamma_k>t|\overline{L}_k=\overline{l}_k,
\overline{A}_{k-1}=\overline{a}_{k-1},T>\tau_k\right)\\
&=&s_{\overline{l}_k,\left(\overline{a}_{k-1},\overline{0}\right)}(t).
\end{eqnarray*}
\end{thm}

\proof
We use backward induction on $k$, starting with $k=K$ and ending with $k=0$.
For $k=K$,
\begin{eqnarray*}
P\bigl(T_K^\gamma>t|\overline{L}_K=\overline{l}_K,
\overline{A}_K=\overline{a}_K,T>\tau_K\bigr)&=&
P\bigl(\gamma_{\overline{l}_K,\overline{a}_K}\left(T\right)>t|
\overline{L}_K=\overline{l}_K,\overline{A}_K=\overline{a}_K,T>\tau_K\bigr)\\
&=&P\bigl(T>\gamma^{-1}_{\overline{l}_K,\overline{a}_K}(t)|
\overline{L}_K=\overline{l}_K,\overline{A}_K=\overline{a}_K,T>\tau_K\bigr)\\
&=&s_{\overline{l}_K,\left(\overline{a}_K,\overline{0}\right)}
\bigl(\gamma^{-1}_{\overline{l}_K,\overline{a}_K}(t)\bigr)\\
&=&s_{\overline{l}_K,\left(\overline{a}_{K-1},\overline{0}\right)}(t).
\end{eqnarray*}
Here the first equality is immediate from the definition of
$T_K^\gamma$, the second follows by the strict monotonicity of the
functions $\gamma$, the third by definition of $s$ and the last by
definition of $\gamma$.

Induction step: we show that if the theorem is true for $k+1$, then it
is also true for $k$. Just as for $k=K$,
\begin{equation*}
P\bigl(T_k^\gamma>t|\overline{L}_k=\overline{l}_k,
\overline{A}_k=\overline{a}_k,T>\tau_k\bigr)
=P\bigl(T_{k+1}^\gamma>\gamma^{-1}_{\overline{l}_k,\overline{a}_k}(t)|
\overline{L}_k=\overline{l}_k,\overline{A}_k=\overline{a}_k,T>\tau_k\bigr).
\end{equation*}
Now we distinguish two possibilities:
$\gamma^{-1}_{\overline{l}_k,\overline{a}_k}(t)\leq\tau_{k+1}$ and
$\gamma^{-1}_{\overline{l}_k,\overline{a}_k}(t)>\tau_{k+1}$. In the
first case, the right side of the preceding display is equal to 
\begin{eqnarray*}
&&P\bigl(T>\gamma^{-1}_{\overline{l}_k,\overline{a}_k}(t)|
\overline{L}_k=\overline{l}_k,\overline{A}_k=\overline{a}_k,T>\tau_k\bigr)\\
&&\hspace{1cm}=s_{\overline{l}_k,\left(\overline{a}_k,\overline{0}\right)}
\bigl(\gamma^{-1}_{\overline{l}_k,\overline{a}_k}(t)\bigr)\\
&&\hspace{1cm}=s_{\overline{l}_k,\left(\overline{a}_{k-1},\overline{0}\right)}(t),
\end{eqnarray*}
where the first equality holds because for
$s\in\left(\tau_k,\tau_{k+1}\right]$ we have that
$\left\{T_{k+1}^\gamma>s\right\}=\left\{T>s\right\}$ by the construction of
$T_{k+1}^\gamma$, and the last equality holds by the definition of $\gamma$.
In the second possibility, i.e.\ if
$\gamma^{-1}_{\overline{l}_k,\overline{a}_k}(t)>\tau_{k+1}$,
\begin{eqnarray*}
\lefteqn{P\bigl(T_{k+1}^\gamma>\gamma^{-1}_{\overline{l}_k,\overline{a}_k}(t)|
\overline{L}_k=\overline{l}_k,\overline{A}_k=\overline{a}_k,T>\tau_k\bigr)}\\
&=&P\bigl(T_{k+1}^\gamma>\tau_{k+1}|\overline{L}_k=\overline{l}_k,
\overline{A}_k=\overline{a}_k,T>\tau_k\bigr)\\
&&\hspace{0.5cm} P\bigl(T_{k+1}^\gamma>\gamma^{-1}_{\overline{l}_k,
\overline{a}_k}(t)|\overline{L}_k=\overline{l}_k,
\overline{A}_k=\overline{a}_k,T>\tau_k,T_{k+1}^\gamma>\tau_{k+1}\bigr)\\
&=&P\left(T>\tau_{k+1}|\overline{L}_k=\overline{l}_k,
\overline{A}_k=\overline{a}_k,T>\tau_k\right)\\
&&\hspace{0.5cm} \sum_{l_{k+1}}\Big\{P\left(L_{k=1}=l_{k+1}|
\overline{L}_{k}=\overline{l}_{k},\overline{A}_k=\overline{a}_k,
T>\tau_{k+1}\right)\\
&&\hspace{0.5cm} \hspace{0.95cm}
P\bigl(T_{k+1}^\gamma>\gamma^{-1}_{\overline{l}_k,\overline{a}_k}(t)|
\overline{L}_{k+1}=\overline{l}_{k+1},\overline{A}_k=\overline{a}_k,
T>\tau_{k+1}\bigr)\Big\}\\
&=&P\left(T>\tau_{k+1}|\overline{L}_k=\overline{l}_k,
\overline{A}_k=\overline{a}_k,T>\tau_k\right)\\
&&\hspace{0.5cm} \sum_{l_{k+1}}\Big\{
P\left(L_{k=1}=l_{k+1}|\overline{L}_{k}=\overline{l}_{k},
\overline{A}_k=\overline{a}_k,T>\tau_{k+1}\right)
s_{\overline{l}_{k+1},\left(\overline{a}_k,\overline{0}\right)}
\bigl(\gamma^{-1}_{\overline{l}_k,\overline{a}_k}(t)\bigr)\Big\}\\
&=&P\left(T>\tau_{k+1}|\overline{L}_k=\overline{l}_k,
\overline{A}_k=\overline{a}_k,T>\tau_k\right)\\
&&\hspace{0.5cm} \sum_{l_{k+1}}\Big\{
P\left(L_{k=1}=l_{k+1}|\overline{L}_{k}=\overline{l}_{k},
\overline{A}_k=\overline{a}_k,T>\tau_{k+1}\right)
s_{\overline{l}_{k+1},\left(\overline{a}_{k-1},\overline{0}\right)}
\left(t\right)\Big\}\\
&=&s_{\overline{l}_{k},\left(\overline{a}_{k-1},\overline{0}\right)}(t),
\end{eqnarray*}
where in the first step we condition on $T_{k+1}^\gamma>\tau_{k+1}$,
in the second we use that
$\left\{T_{k+1}^\gamma>\tau_{k+1}\right\}=\left\{T>\tau_{k+1}\right\}$
and we condition on $L_{k+1}$, the fourth is the induction step, the
fifth follows from the definition of $\gamma$ and the last from
the definition of
$s_{\overline{l}_k,\left(\overline{a}_{k-1},\overline{0}\right)}$.
\Endproof

\end{appendix}

\bigskip\noindent
{\bf Acknowledgement.} This paper is based on an earlier manuscript
by the first author.


\bigskip\noindent
Corresponding author:\\
{\sl Aad van der Vaart}\\
{\sl Department of Mathematics}\\
{\sl Faculty of Sciences}\\
{\sl Vrije Universiteit}\\
{\sl De Boelelaan 1081 a}\\
{\sl 1081 HV Amsterdam}\\
{\sl The Netherlands}\\

\end{document}